\documentclass[11pt,leqno]{article}
\usepackage{a4}
\usepackage[T1]{fontenc}
\usepackage[english]{babel}
\usepackage{latexsym}
\usepackage{amsmath}
\usepackage{amssymb}
\usepackage{mathrsfs}
\usepackage{fig4tex}
\usepackage{color}
\usepackage{cite}
\usepackage{titling}
\usepackage{bbold}
\makeatletter
\@addtoreset{equation}{section} 
\makeatother
\definecolor{move}{rgb}{.3,.1,.8}
\newcommand{\ud}{\mathrm{d}}
\newcommand{\vfi}{\varphi}
\newcommand{\foo}{\vfi_{1,1}}
\newcommand{\fot}{\vfi_{1,2}}
\newcommand{\fto}{\vfi_{2,1}}
\newcommand{\ftt}{\vfi_{2,2}}
\newcommand{\supp}{\mathrm{supp}}
\newcommand{\re}{\mathrm{Re}}
\newcommand{\im}{\mathrm{Im}}
\newcommand{\e}{\mathrm{e}}
\newcommand{\id}{\mathrm{Id}}
\newcommand{\Rn}{\mathbb{R}^{n}}
\newcommand{\R}{\mathbb{R}}
\newcommand{\N}{\mathbb{N}}

\newcommand{\p}{\partial}
\newcommand{\pn}{\partial_{\nu}}
\newcommand{\Ci}{\mathscr{C}^{\infty}}
\newcommand{\Cc}{\mathscr{C}_{c}^{\infty}}

\newenvironment{pr}{\vspace{5pt}\textbf{{\small Proof :}}\\}{\hspace{\stretch{1}}\rule{1ex}{1ex}\vspace{5pt}}
\newtheorem{thm}{Theorem}[section]
\newtheorem{pro}{Proposition}[section]

\newtheorem{lem}{Lemma}[section]
\newtheorem{rem}{Remarks}[section]

\usepackage{fancyhdr}
\fancypagestyle{plain}{\fancyhead{}}
\title{Logarithmic stabilization of the Euler-Bernoulli transmission plate equation with locally distributed Kelvin-Voigt damping}
\author{{FATHI HASSINE}\\ \textit{UR Analysis and Control of PDE (13ES64)}\\ \textit{D\'epartement de Math\'ematiques, Facult\'e des Sciences de Monastir}\\ \textit{Universit\'e de Monastir, 5019 Monastir, Tunisie}\\ \textit{email:} \texttt{fathi.hassine@fsm.rnu.tn}}
\date{}
\begin{document}
\maketitle
\begin{center}
\abstract{In this paper we will study the asymptotic behaviour of the energy decay of a transmission plate equation with locally distributed Kelvin-Voigt feedback. Precisly, we shall prove that the energy decay at least logarithmically over the time. The originality of this method comes from the fact that using a Carleman estimate for a transmission second order system which will be derived from the plate equation to establish a resolvent estimate which provide, by the famous Burq's result~\cite{Bur}, the kind of decay mentionned above.}
\end{center}

\section{Introduction and statement of results}
In recent years, there has been much interest in the stability problems for elastic systems with locally distributed damping. Most of the works were devoted to the viscous damping, i.e., the damping is proportional to the velocity (see for instance~\cite{CFNS} and~\cite{Z}). Structures with local viscoelasticity arise from use of smart material or passive stabilization of structures. However, very little is known about exponential stability for elastic systems with local viscoelastic damping, although there is a fairly deep understanding when the damping is distributed over the entire domain but only for 1-dimension (see~\cite{LL2}). To our knowledge, the first paper in this direction was published in 1998 by Liu and Liu~\cite{LL} where they obtained exponential stability for the Euler-Bernoulli beam equation with local Kelvin-Voigt damping. Noting that in our knowledge there are zero results at least for the multi-dimension Euler-Bernoulli plate equation case. 

Consider a clamped elastic domain in $\Rn$, ($n\geq 2$) which is made of a viscoelastic material with Kelvin-Voigt constitutive relation in which a transmisson effect has been established such a way that the damping is locally effective in only one side the transmission boundary. By the Kirchhoff hypothesis, neglecting the rotatory inertia, the transversal vibration (see~\cite{CLL} for the modeling problem) can be described as follows: Let $\Omega$ and $\Omega_{1}$ be two open, bounded and connected domains with smooth boundary respectively $\Gamma$ and $S$ such that $\Omega_{1}\subset\Omega$ and $\overline{S}\cap\overline{\Gamma}=\emptyset$. We set also $\Omega_{2}=\Omega\backslash\overline{\Omega}_{1}$ which is an open connected domain with boundary $\p\Omega_{2}=\Gamma\cup S$.

We are going to study the following transmission and boundary value problem
\begin{equation}\label{1}
\left\{
\begin{array}{lll}
\p_{t}^{2}u_{1}+\Delta(c_{1}^{2}\Delta u_{1}+a.\Delta\p_{t}u_{1})=0&\textrm{in}&\Omega_{1}\times]0,+\infty[,
\\
\p_{t}^{2}u_{2}+c_{2}^{2}\Delta^{2}u_{2}=0&\textrm{in}&\Omega_{2}\times]0,+\infty[,
\\
u_{1}=u_{2}&\textrm{on}&S\times]0,+\infty[,
\\
\pn u_{1}=\pn u_{2}&\textrm{on}&S\times]0,+\infty[,
\\
c_{1}\Delta u_{1}=c_{2}\Delta u_{2}&\textrm{on}&S\times]0,+\infty[,
\\
c_{1}\pn \Delta u_{1}=c_{2}\pn \Delta u_{2}&\textrm{on}&S\times]0,+\infty[,
\\
u_{1}=0&\textrm{on}&\Gamma\times]0,+\infty[,
\\
\Delta u_{1}=0&\textrm{on}&\Gamma\times]0,+\infty[,
\\
u_{1}(x,0)=u_{1}^{0}(x),\;\p_{t}u_{1}(x,0)=u_{1}^{1}(x)&\textrm{in}&\Omega_{1},
\\
u_{2}(x,0)=u_{2}^{0}(x),\;\p_{t}u_{2}(x,0)=u_{2}^{1}(x)&\textrm{in}&\Omega_{2}.
\end{array}\right.
\end{equation}
Where $\pn$ denotes the unit outward normal vector of $\Omega_{1}$ and $\Omega$ respectively in  $S$ and $\Gamma$, $c_{1}$, $c_{2}$ are strictly positives constants and $a$ is a non negative bounded functions in $\Omega_{1}$ and we suppose that $a$ vanishing near the boundary $S$ such that there exist a non empty open domain $\omega\subset\Omega_{1}$ such that $a$ is strictly positives in $\overline{\omega}$.

The energy of a solution of~\eqref{1} at time $t\geq 0$ is defined by
\begin{equation*}
E(t)=\frac{1}{2}\int_{\Omega_{1}}\Big(|\p_{t}u_{1}(x,t)|^{2}+c_{1}^{2}|\Delta u_{1}(x,t)|^{2}\Big)c_{1}^{-1}\,\ud x+\frac{1}{2}\int_{\Omega_{2}}\Big(|\p_{t}u_{2}(x,t)|^{2}+c_{2}^{2}|\Delta u_{2}(x,t)|^{2}\Big)c_{2}^{-1}\,\ud x.
\end{equation*}
By Green's formula we can prove that for all $\;t_{1},\,t_{2}>0$ we have
\[
E(t_{2})-E(t_{1})=-c_{1}^{-1}\int_{t_{1}}^{t_{2}}\!\!\!\int_{\Omega_{1}}a|\Delta\p_{t}u_{1}(x,t)|^{2}\,\ud x\,\ud t,
\]
and this mean that the energy is decreasing over the time.

We define the operator $\mathcal{A}$ by
$$\mathcal{A}\left(\begin{array}{c}
u_{1}
\\
u_{2}
\\
v_{1}
\\
v_{2}
\end{array}\right)=(v_{1},v_{2},-\Delta(c_{1}^{2}\Delta u_{1}+a\Delta v_{1}),-c_{2}^{2}\Delta^{2} u_{2})$$
in the Hilbert space $\mathcal{H}=X\times H$ where $H=H_{1}\times H_{2}=L^{2}(\Omega_{1},c_{1}^{-1}\,\ud x)\times L^{2}(\Omega_{2},c_{1}^{-1}\,\ud x)$ and
\begin{equation}\label{8}
\begin{split}
X=\big\{(u_{1},u_{2})\in H\,:\,u_{1}\in H^{2}(\Omega_{1}),\,u_{2}\in H^{2}(\Omega_{2}),\,u_{2|\Gamma}=0,\,u_{1\,|S}=u_{2\,|S},\pn u_{1\,|S}=\pn u_{2\,|S}\big\},
\end{split}
\end{equation}
with domain
\begin{align*}
\mathcal{D}(\mathcal{A})=\big\{(u_{1},u_{2},v_{1},v_{2})\in\mathcal{H}\,:\,(v_{1},v_{2},\Delta(c_{1}^{2}\Delta u_{1}+a\Delta v_{1}),c_{2}^{2}\Delta^{2}u_{2})\in\mathcal{H},\,\Delta u_{2\,|\Gamma}=0,
\\
\,c_{1}\Delta u_{1\,|S}=c_{2}\Delta u_{2\,|S},\,c_{1}\pn\Delta u_{1\,|S}=c_{2}\pn\Delta u_{2\,|S}\big\}.
\end{align*}
Now we are able to state our main results
\begin{thm}\label{9}
There exists $C>0$ such that for every $\mu\in \R$ with $|\mu|$ large, we have
\begin{equation}\label{2}
\|(\mathcal{A}-i\mu\,\id)^{-1}\|_{\mathcal{L}(\mathcal{H})}\leq C\e^{C|\mu|}.
\end{equation}
\end{thm}
As an immediate consequence of the previous theorem (see~\cite{Bur} and more recently~\cite{BD}), we get the following rate of decrease of energy
\begin{thm}\label{10}
For any $k\in\N$, there exists a constant $C>0$ such that for any initial data $(u_{1}^{0},u_{2}^{0},u_{1}^{1},u_{2}^{1})\in\mathcal{D}(\mathcal{A}^{k})$, the energy $E(t)$ of the system~\eqref{1} whose solution  $u(x,t)$ is starting from $(u_{1}^{0},u_{2}^{0},u_{1}^{1},u_{2}^{1})$ satisfy
$$E(t)\leq\frac{C}{(\ln(2+t))^{2k}}\|(u_{1}^{0},u_{2}^{0},u_{1}^{1},u_{2}^{1})\|_{\mathcal{D}(\mathcal{A}^{k})}^{2},\quad \forall\; t>0.$$
\end{thm}
\begin{rem}\rm{
\*
\begin{enumerate}
	\item[1)] Under one assumption to the coefficients $c_{1}$ and $c_{2}$, Ammari and Vodev~\cite{AV} have proved an exponential stabilization result for the Euler-Bernoulli transmission plate equation with boundary dissipation.  Again for a transmission model, Ammari and Nicaise~\cite{AN} have proved, under some geometric condition, an exponential stabilization for a coupled damped wave equation with a damped Kirchhoff plate equation.
	\item[2)] To prove Theorem~\ref{9} and Theorem~\ref{10}, we make use the Carleman estimates to obtain information about the resolvent in a boundary domain, the cost is to use phases functions satisfying H\"ormander's assumption. Albano~\cite{A} proved a Carleman estimate for the plate operator, by decomposing the operator as the product of two Schr{\"o}dinger ones and gives for eatch of them the corresponding Carleman estimate then by making together these two estimates we obtain the result. But here we will not need to have a Carleman estimate for the plate equation, namely inspiring from the Albano's decomposition we will derive a second ordre transmission system to which we are going to apply an appropriate Carleman estimates (see section~\ref{b1}) for a suitable phases functions, thus we will obtain the resolvent estimate of Theorem~\ref{9}.
	\item[3)] Theorem~\ref{9} and Theorem~\ref{10} are analogous to those of Fathallah~\cite{I}, in the case of hyperbolic-parabolic coupled system, and Lebeau and Robbiano~\cite{LR} resuts, in the case of scalar wave equation without transmission, but our method is different from their because it consist to use the Carleman estimates directly for the stationary operator without going through the interpolation inequality.
	\item[4)] For various purposes, several authors have focused to the transmission problems where they needed to find a Carleman estimates near the interface, such as the works of Bellassoued~\cite{B} and Fathallah~\cite{I} for the stabilization problems and that also of Le Rousseau and Robbiano~\cite{RR} for a control problem.
	\item[5)] Note that in the case where it has no transmission of the problem~\eqref{1}, Theorem~\ref{9} and Theorem~\ref{10} remain valid and in this case we need only the classical Carleman estimates (see~\cite{LR} and~\cite{LR2}).
\end{enumerate}
}
\end{rem}

In this paper $C$ will always be a generic positive constant whose value may be different from one line to another.

The outline of this paper is as follow. In section~\ref{11} we prove the well-Posedness of the problem~\eqref{1}, in section~\ref{b1} we give a global Carleman estimate and we will constract a suitable phases functions and in section~\ref{d1} we prove the resolvent estimate gived by Theorem~\ref{9}.
\section{Well-Posedness of the problem}\label{11}
To prove the Well-Posedness of the problem~\eqref{1} we are going to use the semigroups theory. Our strategy consiste to write the equations as a Cauchy problem with an operator which generates a semigroup of contractions.

Throughout this paper, we denote the inner product in the space $H=H_{1}\times H_{2}$ by
\begin{equation*}
\left\langle\left(\begin{array}{c}
u_{1}
\\
u_{2}
\end{array}\right),\left(\begin{array}{c}
v_{1}
\\
v_{2}
\end{array}\right)\right\rangle_{H}=\int_{\Omega_{1}}u_{1}(x)\overline{v_{1}(x)}c_{1}^{-1}\,\ud x+\int_{\Omega_{2}}u_{2}(x)\overline{v_{2}(x)}c_{2}^{-1}\,\ud x,
\end{equation*}
The Cauchy problem is written in the following form
$$\left\{\begin{array}{ll}
\p_{t}\left(\begin{array}{c}
u_{1}
\\
u_{2}
\\
v_{1}
\\
v_{2}
\end{array}\right)(t)=\mathcal{A}\left(\begin{array}{c}
u_{1}
\\
u_{2}
\\
v_{1}
\\
v_{2}
\end{array}\right)(t)&t\in]0,+\infty[,
\\
\left(\begin{array}{l}
u_{1}
\\
u_{2}
\\
v_{1}
\\
v_{2}
\end{array}\right)(0)=\left(\begin{array}{l}
u_{1}^{0}
\\
u_{2}^{0}
\\
u_{1}^{1}
\\
u_{2}^{1}
\end{array}\right).&
\end{array}\right.$$

Now we have to specify the functional space and the domain of the operator $\mathcal{A}$. In the space $H$ we define the operator $G$ by
\[
G\left(\begin{array}{l}
u_{1}
\\
u_{2}\end{array}\right)=(-c_{1}\Delta u_{1},-c_{2}\Delta u_{2})\qquad \forall\;(u_{1},u_{2})\in\mathcal{D}(G)
\]
with domain $\mathcal{D}(G)=X$ defined in~\eqref{8}. The space $X$ is equipped with the norm
$$\|(u_{1},u_{2})\|_{X}=\|G(u_{1},u_{2})\|_{H}$$
and we defined the graph norm of $G$ by
$$\|(u_{1},u_{2})\|_{gr(G)}^{2}=\|(u_{1},u_{2})\|_{H}^{2}+\|G(u_{1},u_{2})\|_{H}^{2}$$
then we have the following
\begin{pro}
$(X,\|\,.\,\|_{X})$ is a Hilbert space with a norm equivalent to the graph norm of $G$.
\end{pro}
\begin{pr}
It is well known that if $G$ is a colsed operator then $(X,\|\,.\,\|_{gr(G)})$ is a Hilbert space. Thus to prove the proposition it suffices to show that $G$ is closed and both norms are equivalent.
\\
By Green's formula and Poincar\'e inequality it is easy to show that there exists $C>0$ such that
$$\left\langle G\left(\begin{array}{l}
u_{1}
\\
u_{2}\end{array}\right),\left(\begin{array}{l}
u_{1}
\\
u_{2}\end{array}\right)\right\rangle_{H}=\|\nabla u_{1}\|_{L^{2}(\Omega_{1})}^{2}+\|\nabla u_{2}\|_{L^{2}(\Omega_{2})}^{2}\geq C\|(u_{1},u_{2})\|_{H}^{2}\quad\forall\,(u_{1},u_{2})\in X.$$
Then $G$ is a strictly positive operator and we have
$$\|G(u_{1},u_{2})\|_{H}.\|(u_{1},u_{2})\|_{H}\geq\left\langle G\left(\begin{array}{l}
u_{1}
\\
u_{2}\end{array}\right),\left(\begin{array}{l}
u_{1}
\\
u_{2}\end{array}\right)\right\rangle_{H}\geq C\|(u_{1},u_{2})\|_{H}^{2}\quad\forall\,(u_{1},u_{2})\in X$$
which prove the equivalence between the two norms.
\\
Now since $G$ is positive then by in~\cite[Propsition 3.3.5]{TW}, $-G$ is m-dissipative and thus $G$ is a closed operator. This completes the proof.
\end{pr}

This last result allows us to properly define the functional space of the operator $\mathcal{A}$.
\begin{pro}\label{7}
The two spaces $(X,\|\,.\,\|_{2})$ and $(X,\|\,.\,\|_{X})$ are algebraically and topologically the same. Where we have defined $\|\,.\,\|_{2}$ by
$$\|(u_{1},u_{2})\|_{2}^{2}=\|u_{1}\|_{H^2(\Omega_{1})}^{2}+\|u_{2}\|_{^{H^2(\Omega_{2})}}^{2},\qquad\forall\,(u_{1},u_{2})\in X.$$
\end{pro}
\begin{pr}
We have only to prove that the two norms are equivalent.
\\
First, we note that $(X,\|\,.\,\|_{2})$ is a Hilbert space because $X$ is a closed subspace of $H^{2}(\Omega_{1}\cup\Omega_{2})$, in addition we have
$$\|(u_{1},u_{2})\|_{X}^{2}=\|\Delta u_{1}\|_{L^{2}(\Omega_{1})}^{2}+\|\Delta u_{2})\|_{L^{2}(\Omega_{2})}^{2}\leq C\|(u_{1},u_{2})\|_{2}^{2}\quad\forall\,u\in X,$$
and while $(X,\|\,.\,\|_{X})$ is also a Hilbert space, then according to the Banach theorem (see~\cite[Corollary 9.2.3]{YVA}) the two norms are equivalent.
\end{pr}

We set $\mathcal{H}=X\times H$ the Hilbert space with the norm
\[
\|(u_{1},u_{2},v_{1},v_{2})\|^{2}=\|(u_{1},u_{2})\|_{X}^{2}+\|(v_{1},v_{2})\|_{H}^{2}\qquad\forall\, (u_{1},u_{2},v_{1},v_{2})\in\mathcal{H},
\]
and we recall that the domain of the operator $\mathcal{A}$ is defined by
\begin{align*}
\mathcal{D}(\mathcal{A})=\big\{(u_{1},u_{2},v_{1},v_{2})\in\mathcal{H}\,:\,(v_{1},v_{2},\Delta(c_{1}^{2}\Delta u_{1}+a\Delta v_{1}),c_{2}^{2}\Delta^{2}u_{2})\in\mathcal{H},\,\Delta u_{2\,|\Gamma}=0,
\\
\,c_{1}\Delta u_{1\,|S}=c_{2}\Delta u_{2\,|S},\,c_{1}\pn\Delta u_{1\,|S}=c_{2}\pn\Delta u_{2\,|S}\big\}.
\end{align*}
\begin{thm}
Under the above assumptions, the operator $\mathcal{A}$ is m-dissipative and especially it generates a strongly semigroup of contractions in $\mathcal{H}$.
\end{thm}
\begin{pr}
According to Lumer-Phillips theorem (see for exemple~\cite[p.103]{TW}) we have only to prove that $\mathcal{A}$ is m-dissipative.
\\
Let $(u_{1},u_{2},v_{1},v_{2})\in\mathcal{D}(\mathcal{A})$ then by Green's formula we have
\begin{eqnarray*}
\re\left\langle\mathcal{A}\left(\begin{array}{l}
u_{1}
\\
u_{2}
\\
v_{1}
\\
v_{2}
\end{array}\right),\left(\begin{array}{l}
u_{1}
\\
u_{2}
\\
v_{1}
\\
v_{2}
\end{array}\right)\right\rangle_{\mathcal{H}}&=&\re\left\langle\left(\begin{array}{c}
v_{1}
\\
v_{2}
\\
-\Delta(c_{1}^{2}\Delta u_{1}+a\Delta v_{1})
\\
-c_{2}^{2}\Delta^{2}u_{2}
\end{array}\right),\left(\begin{array}{l}
u_{1}
\\
u_{2}
\\
v_{1}
\\
v_{2}
\end{array}\right)\right\rangle_{\mathcal{H}}
\\
&=&-c_{1}\|a^{\frac{1}{2}}\Delta v_{1}\|_{L^{2}(\Omega_{1})}^{2}\leq 0.
\end{eqnarray*}
This shows that $\mathcal{A}$ is dissipative.
\\
Let now $(f_{1},f_{2},g_{1},g_{2})\in\mathcal{H}$ and our purpose is to find a couple $(u_{1},u_{2},v_{1},v_{2})\in\mathcal{D}(\mathcal{A})$ such that
$$\left(\id-\mathcal{A}\right)\left(\begin{array}{l}
u_{1}
\\
u_{2}
\\
v_{1}
\\
v_{2}
\end{array}\right)=\left(\begin{array}{c}
u_{1}-v_{1}
\\
u_{2}-v_{2}
\\
v_{1}+\Delta(c_{1}^{2}\Delta u_{1}+a\Delta v_{1})
\\
v_{2}+c_{2}^{2}\Delta^{2}u_{2}
\end{array}\right)=\left(\begin{array}{l}
f_{1}
\\
f_{2}
\\
g_{1}
\\
g_{2}
\end{array}\right)$$
more explicitly we have to find $(u_{1},u_{2},v_{1},v_{2})\in\mathcal{D}(\mathcal{A})$ such that
\begin{equation*}
\left\{\begin{array}{l}
v_{1}=u_{1}-f_{1}
\\
v_{2}=u_{2}-f_{2}
\\
u_{1}+\Delta((c_{1}^{2}+a)\Delta u_{1}-a\Delta f_{1})=f_{1}+g_{1}
\\
u_{2}+c_{2}^{2}\Delta^{2} u_{2}=f_{2}+g_{2}.
\end{array}\right.
\end{equation*}
First note that, by Riesz representation theorem, there exists a unique $(u_{1},u_{2})\in X=\mathcal{D}(G)$ such that for all $(\vfi_{1},\vfi_{2})\in X$ we have
\begin{equation}\label{4}
\begin{split}
\langle f_{1}+g_{1},\vfi_{1}\rangle_{L^{2}(\Omega_{1})}+\langle f_{2}+g_{2},\vfi_{2}\rangle_{L^{2}(\Omega_{2})}+\langle a\Delta f_{1},\Delta\vfi_{1}\rangle_{L^{2}(\Omega_{1})}=\langle u_{1},\vfi_{1}\rangle_{L^{2}(\Omega_{1})}
\\
+\langle u_{2},\vfi_{2}\rangle_{L^{2}(\Omega_{2})}+\langle(c_{1}^{2}+a)\Delta u_{1},\Delta\vfi_{1}\rangle_{L^{2}(\Omega_{1})}+c_{2}^{2}\langle\Delta u_{2},\Delta\vfi_{2}\rangle_{L^{2}(\Omega_{2})}.
\end{split}
\end{equation}
In particular for all $(\vfi_{1},\vfi_{2})\in\Cc(\Omega_{1})\times\Cc(\Omega_{2})$ the expression~\eqref{4} yields
\begin{eqnarray*}
\langle\Delta((c_{1}^{2}+a)\Delta u_{1}-a\Delta f_{1})+(u_{1}-f_{1}-g_{1}),\vfi_{1}\rangle_{D'(\Omega_{1})}=0,
\\
\langle c_{2}^{2}\Delta^{2}u_{2}+(u_{2}-f_{2}-g_{2}),\vfi_{2}\rangle_{D'(\Omega_{2})}=0.
\end{eqnarray*}
then we obtain
\begin{equation}\label{5}
\begin{split}
u_{1}+\Delta((c_{1}^{2}+a)\Delta u_{1}-a\Delta f_{1})=f_{1}+g_{1}\quad\text{in}\;L^{2}(\Omega_{1}),
\\
u_{2}+c_{2}^{2}\Delta^{2} u_{2}=f_{2}+g_{2}\quad\text{in}\;L^{2}(\Omega_{2}).
\end{split}
\end{equation}
Now if we return again to the expression~\eqref{4} then by Green's formula we can write it as follows
\begin{align*}
\langle\Delta((c_{1}^{2}+a)\Delta u_{1}-a\Delta f_{1})+(u_{1}-f_{1}-g_{1}),\vfi_{1}\rangle_{L^{2}(\Omega_{1})}+
\langle c_{2}^{2}\Delta^{2}u_{2}+(u_{2}-f_{2}-g_{2}),\vfi_{2}\rangle_{L^{2}(\Omega_{2})}
\\
=-\langle c_{2}\Delta u_{2},\pn\vfi_{2}\rangle_{L^{2}(\Gamma)}-\langle c_{1}\Delta u_{1}\pn\vfi_{1}\rangle_{L^{2}(S)}+\langle c_{2}\Delta u_{2},\pn\vfi_{1}\rangle_{L^{2}(S)}
\\
+\langle c_{1}\pn\Delta u_{1},\vfi_{2}\rangle_{L^{2}(S)}-\langle c_{2}\pn\Delta u_{2},\vfi_{2}\rangle_{L^{2}(S)}.
\end{align*}
then by~\eqref{5} we get for all $(\vfi_{1},\vfi_{2})\in X$ that
\begin{equation*}
\langle c_{1}\pn\Delta u_{1}-c_{2}\pn\Delta u_{2},\vfi_{1}\rangle_{L^{2}(S)}-\langle c_{1}\Delta u_{1}-c_{2}\Delta u_{2},\pn\vfi_{1}\rangle_{L^{2}(S)}-\langle c_{2}\Delta u_{2},\pn\vfi_{2}\rangle_{L^{2}(\Gamma)}=0,
\end{equation*}
which yields the following equalities
$$c_{1}\Delta u_{1\,|S}=c_{2}\Delta u_{2\,|S},\; c_{1}\pn \Delta u_{1\,|S}=c_{2}\pn \Delta u_{2\,|S},\; \Delta u_{2\,|\Gamma}=0.$$
And this concludes the proof.
\end{pr}

One consequence of this last result is that if we assume that $(u_{1}^{0},u_{2}^{0},u_{1}^{1},u_{2}^{1})\in\mathcal{D}(\mathcal{A})$, there exists a unique solution of~\eqref{1} which can be expressed by means of a semigroup on $\mathcal{H}$ as follows
\begin{equation}\label{6}
\left(\begin{array}{c}
u_{1}
\\
u_{2}
\\
\p_{t}u_{1}
\\
\p_{t}u_{2}
\end{array}\right)=e^{t\mathcal{A}}\left(\begin{array}{c}
u_{1}^{0}
\\
u_{2}^{0}
\\
u_{1}^{1}
\\
u_{2}^{1}
\end{array}\right)
\end{equation}
where $e^{t\mathcal{A}}$ is the $C_{0}$-semigroup of contractions generates by the operator $\mathcal{A}$. And we have the following regularity of the solution
$$
\left(\begin{array}{c}
u_{1}
\\
u_{2}
\\
\p_{t}u_{1}
\\
\p_{t}u_{2}
\end{array}\right)
\in C([0,+\infty[,\mathcal{D}(\mathcal{A}))\cap C^{1}([0,+\infty[,\mathcal{H}).$$
\\
And if $(u_{1}^{0},u_{2}^{0},u_{1}^{1},u_{2}^{1})\in\mathcal{H}$, the function $(u_{1}(t),u_{2}(t))$ given by~\eqref{6} is the mild solution of~\eqref{1} and it lives in $C([0,+\infty[,\mathcal{H})$.
\section{Carleman estimate and construction weight functions}\label{b1}
\subsection{Carleman estimate}
We consider tow open and disjoint domains $\mathcal{O}_{1}$ and $\mathcal{O}_{2}$ in which we define respectively the second order elliptic semi-classical operators $P_{1}=-h^{2}\Delta-\alpha_{1}h$ and $P_{2}=-h^{2}\Delta-\alpha_{2}h$ with principal symbol $p(x,\xi)=|\xi|^{2}$ where $h$ is a very small semi-classical parmeter and $\alpha_{1},\,\alpha_{2}\in\R$, and we suppose that $\p\mathcal{O}_{1}=\gamma\cup\gamma_{1}$, $\p\mathcal{O}_{2}=\gamma\cup\gamma_{2}$ and $\overline{\gamma}_{1}\cap\overline{\gamma}_{0}=\overline{\gamma}_{2}\cap\overline{\gamma}_{0}=\emptyset$.

Let $\vfi_{1}\in\Ci(\overline{\mathcal{O}}_{1})$ and $\vfi_{2}\in\Ci(\overline{\mathcal{O}}_{2})$ tow real value functions. We define the two adjoint operators $P_{\vfi_{1}}=\e^{\vfi_{1}/h}P_{1}\e^{\vfi_{2}/h}$ and $P_{\vfi_{1}}=\e^{\vfi_{1}/h}P_{1}\e^{\vfi_{2}/h}$ of principal symbol respectively $p_{1}(x,\xi)=p(x,\xi+i\nabla\vfi_{1})$ and $p_{2}(x,\xi)=p(x,\xi+i\nabla\vfi_{2})$.

By denoting $\pn$ the unit outward normal vector of $\mathcal{O}_{1}$ and $\mathcal{O}_{2}$ respectively in  $\gamma\cup\gamma_{1}$ and $\gamma_{2}$ we assume that the weight function $\vfi_{1}$ and $\vfi_{2}$ satisfies
\begin{enumerate}
	\item[1)] $|\nabla\vfi_{1}|(x)>0,\;\forall\,x\in\overline{\mathcal{O}}_{1}$ and $|\nabla\vfi_{2}|(x)>0,\;\forall\,x\in\overline{\mathcal{O}}_{2}$,
	\item[2)] $\pn\vfi_{1\,|\gamma_{1}}\neq 0$ and $\pn\vfi_{2\,|\gamma_{2}}<0$,
	\item[3)] $\vfi_{1\,|\gamma}=\vfi_{2\,|\gamma}$,
	\item[4)]$(\pn\vfi_{1})_{|\gamma}< 0$, $(\pn\vfi_{2})_{|\gamma}<0$ and $(\pn\vfi_{1})_{|\gamma}^{2}-(\pn\vfi_{2})_{|\gamma}^{2}>0$,
	\item[5)]The sub-ellipticity condition respectively in $\overline{\mathcal{O}}_{1}$ and $\overline{\mathcal{O}}_{2}$
	\begin{equation*}
	\begin{array}{l}
	\forall\,(x,\xi)\in\overline{\mathcal{O}}_{1}\times\Rn;\; p_{\vfi_{1}}(x,\xi)=0\Longrightarrow\{\re(p_{\vfi_{1}}),\im(p_{\vfi_{1}})\}(x,\xi)>0,
	\\
	\forall\,(x,\xi)\in\overline{\mathcal{O}}_{2}\times\Rn;\; p_{\vfi_{2}}(x,\xi)=0\Longrightarrow\{\re(p_{\vfi_{2}}),\im(p_{\vfi_{2}})\}(x,\xi)>0.
	\end{array}
	\end{equation*}
\end{enumerate}
The Carleman estimate corresponding to the following transmission boundary value problem
\begin{equation}\label{b2}
\left\{\begin{array}{ll}
\displaystyle-\Delta w_{1}-\frac{\alpha_{1}}{h}w_{1}=f_{1}&\text{in  }\mathcal{O}_{1}
\\
\\
\displaystyle-\Delta w_{2}-\frac{\alpha_{2}}{h}w_{2}=f_{2}&\text{in  }\mathcal{O}_{2}
\\
w_{1}=w_{2}+e_{1}&\text{on  }\gamma
\\
\pn w_{1}=\pn w_{2}+e_{2}&\text{on  }\gamma
\\
w_{2}=0&\text{on  }\gamma_{2}
\end{array}\right.
\end{equation}
is gived in the following
\begin{thm}~\cite[Theorem 2.1]{RR}\label{b3}
Under the above assumptions on the weight functions $\vfi_{1}$ and $\vfi_{2}$, there exists $h_{0}>0$ and $C>0$ such that
\begin{equation}\label{b4}
\begin{split}
&h\|\e^{\vfi_{1}/h}w_{1}\|_{L^{2}(\mathcal{O}_{1})}^{2}+h^{3}\|\e^{\vfi_{1}/h}\nabla w_{1}\|_{L^{2}(\mathcal{O}_{1})}^{2}+h|\e^{\vfi_{1}/h}w_{1}|_{L^{2}(\gamma)}^{2}+h^{3}|\e^{\vfi_{1}/h}\nabla w_{1}|_{L^{2}(\gamma)}^{2}+
\\
&h^{3}|\e^{\vfi_{1}/h}\pn w_{1}|_{L^{2}(\gamma)}^{2}+h\|\e^{\vfi_{2}/h}w_{2}\|_{L^{2}(\mathcal{O}_{2})}^{2}+h^{3}\|\e^{\vfi_{2}/h}\nabla w_{2}\|_{L^{2}(\mathcal{O}_{2})}^{2}+h|\e^{\vfi_{2}/h}w_{2}|_{L^{2}(\gamma)}^{2}+
\\
&h^{3}|\e^{\vfi_{2}/h}\nabla w_{2}|_{L^{2}(\gamma)}^{2}+h^{3}|\e^{\vfi_{2}/h}\pn w_{2}|_{L^{2}(\gamma)}^{2}\leq C(h^{4}\|\e^{\vfi_{1}/h}f_{1}\|_{L^{2}(\mathcal{O}_{1})}^{2}+h^{4}\|\e^{\vfi_{2}/h}f_{2}\|_{L^{2}(\mathcal{O}_{2})}^{2}+
\\
&h|\e^{\vfi_{1}/h}w_{1}|_{L^{2}(\gamma_{1})}^{2}+h^{3}|\e^{\vfi_{2}/h}\pn w_{1}|_{L^{2}(\gamma_{1})}^{2}+h|\e^{\vfi_{1}/h} e_{1}|_{L^{2}(\gamma)}^{2}+h^{3}|\e^{\vfi_{1}/h}\nabla e_{1}|_{L^{2}(\gamma)}^{2}+h^{3}|\e^{\vfi_{1}/h}e_{2}|_{L^{2}(\gamma)}^{2})
\end{split}
\end{equation}
for all $w_{1}\in\Ci(\overline{O}_{1})$ and $w_{2}\in\Ci(\overline{O}_{2})$ satisfing the system~\eqref{b2} and $h\in]0,h_{0}]$.
\end{thm}
\begin{rem}\rm{
\*
\begin{enumerate}
	\item[1)] If the function $w_{1}$ is supported away from $\gamma_{1}$ the estimate~\eqref{b4} is allows true even if we don't assume that $(\pn\vfi_{1})_{|\gamma_{1}}\neq 0$, while the proof of Theorem~\ref{b3} is local.
	\item[2)] We can not assume that $(\pn\vfi_{1})_{|\gamma_{1}}<0$ (it means $\pn\vfi_{1}<0$ in whole $\p\mathcal{O}_{1}$), otherwise the weight function attain his global maximum in $\mathcal{O}_{1}$ and thus our srtategy of the construction of the phases is fails (see below).
\end{enumerate}
}
\end{rem}
\subsection{Weight function's construction}
\begin{figure}[htbp]
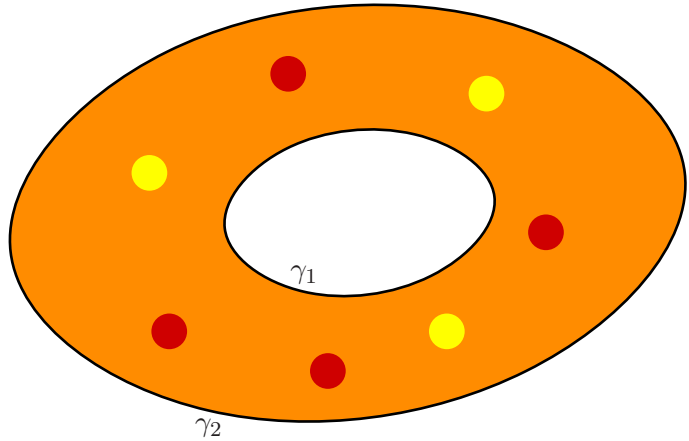

\figinit{1.5pt}
\figpt 0:(20,0)
\figpt 1:(-70,-10)
\figpt 2:(-20,-50)
\figpt 3:(60,-40)
\figpt 4:(100,10)
\figptsym 5:=3/1,4/
\figptsym 6:=2/1,4/
\figpthom 14:=0/1,0.6/
\figpthom 15:=0/2,0.6/
\figpthom 16:=0/3,0.6/
\figpthom 17:=0/4,0.6/
\figpthom 18:=0/5,0.6/
\figpthom 19:=0/6,0.6/
\figpt 7:(50,30)
\figpt 8:(0,35)
\figpt 9:(-35,10)
\figpt 10:(-30,-30)
\figpt 11:(65,-5)
\figpt 12:(10,-40)
\figpt 13:(40,-30)
\psbeginfig{}
\psset(width=1)
\psset(fillmode=yes,color=\DarkOrangergb)
\pscurve[1,2,3,4,5,6,1,2,3]
\psset(fillmode=yes,color=0.6)
\psset(fillmode=yes,color=1)
\pscurve[14,15,16,17,18,19,14,15,16]
\psset(color=1)
\pscirc 7(4)
\pscirc 8(4)
\pscirc 9(4)
\pscirc 10(4)
\pscirc 11(4)
\pscirc 12(4)
\pscirc 13(4)
\psset(fillmode=no)
\psset(color=\Blackrgb)
\pscurve[1,2,3,4,5,6,1,2,3]
\pscurve [14,15,16,17,18,19,14,15,16]
\psset(fillmode=yes,color=0.8\Redrgb)
\pscirc 8(4.5)
\pscirc 10(4.5)
\pscirc 11(4.5)
\pscirc 12(4.5)
\psset(color=\Yellowrgb)
\pscirc 7(4.5)
\pscirc 9(4.5)
\pscirc 13(4.5)
\psendfig
\figvisu{\figBoxA}{}{
\figwriten 15:$\gamma_{1}$(2)
\figwrites 2:$\gamma_{2}$(2)
}
\centerline{\box\figBoxA}
\centerline{\box\figBoxA}
\caption{The domains of the weight functions $\varphi_{1}$ and $\psi_{1}$ (in yellow and orange), $\varphi_{2}$ and $\psi_{2}$ (in red and orange) where they have not critical points.}
\label{fig2}
\end{figure}

In this section we will try to find two phases that satisfies the H{\"o}rmander's condition except in a finite number of ball where one of them do not satisfies this condition the second does and is strictly greater. Note that this result is similar to the Burq's one~\cite[Proposition 3.2]{Bur}, but here we give a new proof due to F.~Laudenbach. Then we will adapte this result to our case to constructe a suitable weight functions that will be needed in the following section. The main ingredient of this section is the following one.
\begin{pro}\label{b5}
Let $\mathcal{O}$ be a bounded open subset with boundary $\gamma=\gamma_{1}\cup\gamma_{2}$ where $\overline{\gamma}_{1}\cap\overline{\gamma}_{2}=\emptyset$ , then there exists two real functions $\psi_{1},\,\psi_{2}\in\Ci(\mathcal{O})$ and continous on $\overline{\mathcal{O}}$ satisfying for $k=1,2$ that $(\pn\psi_{k})_{|\gamma_{1}}<0$ and $(\pn\psi_{k})_{|\gamma_{2}}>0$  having only degenerate critical points (of finite number) such that when $\nabla\psi_{k}=0$ then $\nabla\psi_{\sigma(k)}\neq 0$ and $\psi_{\sigma(k)}>\psi_{k}$. Where $\sigma$ is the permutation of the set $\{1,2\}$ different from the identity. 
\end{pro}
\begin{rem}\label{b6}
\rm{
\*
\begin{enumerate}
	\item[1)] One consequence of Proposition~\ref{b3} is that for $k=1,2$ we can find a finite number of points $x_{kj_{k}}$ and $j_{k}=1,\ldots,N_{k}$ and $\epsilon>0$ such that $B(x_{kj_{k}},2\epsilon)\subset\overline{\mathcal{O}}$ and $B(x_{1j_{1}},2\epsilon)\cap B(x_{2j_{2}},2\epsilon)=\emptyset$, for all $k=1,2$ and $j_{k}=1,\ldots,N_{k}$ and in $B(x_{kj_{k}},2\epsilon)$ we have $\psi_{\sigma(k)}>\psi_{k}$ (See Figure~\ref{fig2}).
	\item[2)] For $\lambda>0$ large enough the weight functions $\vfi_{k}=\e^{\lambda\psi_{k}}$ satisfy the H{\"o}rmander's condition in $\displaystyle U_{k}=\mathcal{O}\bigcap\left(\bigcup_{j_{k}=1}^{N_{k}}B(x_{kj_{k}},\epsilon)\right)^{c}$. Indeed, we have only to prove that for an open bounded subset $U\in\Rn$ and if $\psi\in\Ci(\overline{U})$ satisfying $|\nabla\psi|\geq C$ in $\overline{U}$ and $\vfi=\e^{\lambda\psi}$ we have $\{\re(p_{\vfi}),\im(p_{\vfi})\}(x,\xi)\geq C'$ in $\overline{U}\times\Rn$ for $\lambda>0$ large enough. We have
$$
\left\{\begin{array}{c}
\nabla\vfi=\lambda\e^{\lambda\psi}\nabla\psi\;\text{ and }\;\vfi''=\e^{\lambda\psi}(\lambda\nabla\psi.{}^{t}\nabla\psi+\lambda\psi'')
\\
p_{\vfi}(x,\xi)=0\Longrightarrow\langle\xi,\nabla\vfi\rangle=0\text{ and }|\xi|^{2}=|\nabla\vfi|^{2}
\end{array}\right.
$$
then we obtain
\begin{eqnarray*}
\{\re(p_{\vfi}),\im(p_{\vfi})\}(x,\xi)&=&4\lambda\e^{\lambda\psi}\,{}^{t}\xi.\psi''.\xi+4\e^{3\lambda\psi}(\lambda^{4}|\nabla\psi|^{2}+\lambda^{3}\,{}^{t}\nabla\psi.\psi''.\nabla\psi)
\\
&=&4\e^{3\lambda\psi}(\lambda^{4}|\nabla\psi|^{2}+O(\lambda^{3})).
\end{eqnarray*}
Which conclude the result.
	\item[3)] In general, Proposition~\ref{b3} is also true for any smooth manifold with boundary which the latter is the disjoint union of two open and closed submanifolds.
\end{enumerate}
}
\end{rem}
\begin{pr}
While the Morse functions are dense (for the $\Ci$ topology) in the set of $\Ci$ functions then we can find $\psi_{1}$ a Morse function such that $(\pn\psi_{1})_{|\gamma_{1}}<0$ and $(\pn\psi_{1})_{|\gamma_{2}}>0$. We can suppose that $\psi_{1}$ have no local maximum in $\mathcal{O}$ (The proceeding of the elimination of the maximum is described by Burq~\cite[Appendix A]{Bur}, we can see also~\cite[Theorem 8.1]{Mi} and~\cite[Lemma 2.6]{L}).

Let $c$ be a critical point of $\psi_{1}$ while its index is different from $n$ then we can find a $\Ci$ arc $\gamma_{c}:[-1,1]\rightarrow\Omega$ such that $\gamma_{c}(0)=c$ and $\psi_{1}(\gamma_{c}(1))=\psi_{1}(\gamma_{c}(-1))>\psi_{1}(c)$. We do this construction for all the critical points of $\psi_{1}$ so that all the arcs are mutually disjoint. Hence, this allows us to find a vector field $X$ in $\mathcal{O}$, vanishing near the boundary of $\mathcal{O}$ such that for all critical points $c$ of $\psi_{1}$ we have
$$X(\gamma_{c}(t))=\stackrel{.}{\gamma}_{c}(t),$$
where $\stackrel{.}{\gamma}$ stand for the time derivative.

We denote $\phi_{t}$ its flow:
$$\stackrel{.}{\phi_{t}}(x)=X(\phi_{t}(x)),$$
and we set $\psi_{2}=\psi_{1}\circ\phi_{1}$, thus $\psi_{1}$ and $\psi_{2}$ satisfy the required properties. Indeed, since $X\equiv0$ near the boundary $\gamma_{1}$ and $\gamma_{2}$ which mean that $\phi_{t}(x)=x$ near $\gamma_{1}$ and $\gamma_{2}$ then $\pn\psi_{1\,|\gamma_{1}}=\pn\psi_{2\,|\gamma_{1}}$ and $\pn\psi_{1\,|\gamma_{2}}=\pn\psi_{2\,|\gamma_{2}}$. If $c$ is a critical point of $\psi_{1}$ then we have $\psi_{2}(c)=\psi_{1}(\gamma_{c}(1))>\psi_{1}(c)$, and if $c'$ is a critical point of $\psi_{2}$ then $c'=\phi_{-1}(c)$ where $c$ is a critical point of $\psi_{1}$ and we have $\psi_{2}(c')=\psi_{1}(\phi_{1}\circ\phi_{-1}(c))=\psi_{1}(c)<\psi_{1}(\phi_{-1}(c))=\psi_{1}(c')$ by the construction of $\gamma_{c}$.
\end{pr}

Now if we return to our geometric baseline as described in the introduction of this paper then according to Proposition~\ref{b5} and Remark~\ref{b6} and by noting $\widetilde{\Omega}_{1}=\Omega_{1}\backslash \overline{B}_{r}$ where $B_{r}$ is an open ball of $\Omega_{1}$ with radius $r>0$ such that $\overline{B}_{r}\subset\Omega_{1}$ we can find four phases $\foo$, $\fot$, $\fto$ and $\ftt$ verifying the H\"ormander's condition respectively in $\displaystyle U_{1,1}=\widetilde{\Omega}_{1}\bigcap\left(\bigcup_{j=1}^{N_{11}}B(x_{11}^{j},\epsilon)\right)^{c}$, $\displaystyle U_{1,2}=\widetilde{\Omega}_{1}\bigcap\left(\bigcup_{j_{2}=1}^{N_{12}}B(x_{12}^{j},\epsilon)\right)^{c}$, $\displaystyle U_{2,1}=\Omega_{2}\bigcap\left(\bigcup_{j_{1}=1}^{N_{21}}B(x_{21}^{j},\epsilon)\right)^{c}$ and $\displaystyle U_{2,2}=\Omega_{2}\bigcap\left(\bigcup_{j_{2}=1}^{N_{22}}B(x_{22}^{j},\epsilon)\right)^{c}$ such that $|\nabla\foo|>0$ in $U_{1,1}$, $|\nabla\fot|>0$ in $U_{1,2}$, $|\nabla\fto|>0$ in $U_{2,1}$ and $|\nabla\ftt|>0$ in $U_{2,2}$, moreover $\vfi_{1,k}<\vfi_{1,\sigma(k)}$ in $B(x_{1k}^{j},2\epsilon)$ for all $j=1,\ldots,N_{1,k}$ and $\vfi_{2,k}<\vfi_{2,\sigma(k)}$ in $B(x_{2k}^{j},2\epsilon)$ for all $j=1,\ldots,N_{2,k}$. Furthermore we have also for all $k=1,2$
$$
(\pn\vfi_{1,k})_{|S}<0,\quad(\pn\vfi_{2,k})_{|S}<0 \textrm{ and }(\pn\vfi_{2,k})_{|\Gamma}<0.
$$
We can suppose also that $\vfi_{1,k\,|S}=\vfi_{2,k\,|S}$, and by argument of density we can suppose also that
$$
(\pn\vfi_{1,k})_{|S}^{2}-(\pn\vfi_{2,k})_{\,|S}^{2}>0.
$$
And this concludes the construction of weight functions that will be used in next section.
\section{Resolvent estimate}\label{d1}
The purpose of this section is to find an estimate of the resolvent $(\mathcal{A}-i\mu\id)^{-1}$ where $\mu$ is a real number such that $|\mu|$ is large enough. More precisely we prove that $\|(\mathcal{A}-i\mu\id)^{-1}\|_{\mathscr{L}(\mathcal{H})}\leq C\e^{C|\mu|}$ which imply the weak energy decay of the solution of the equation~\eqref{1}.

The main idea consiste to applying the Carleman estimates for a second order elliptic transmission system which is derived from the plate equation and this is what comes from the originality of our work, it means we prove the stability result for a system of fourth order by using an estimate of Carleman of second order only.

Let $(f_{1},f_{2},g_{1},g_{2})\in\mathcal{H}$ and $(u_{1},u_{2},v_{1},v_{2})\in \mathcal{D}(\mathcal{A})$ such that
$$(\mathcal{A}-i\mu\id)\left(\begin{array}{c}
u_{1}
\\
u_{2}
\\
v_{1}
\\
v_{2}
\end{array}\right)=\left(\begin{array}{c}
f_{1}
\\
f_{2}
\\
g_{1}
\\
g_{2}
\end{array}\right),$$
then we get the following boundary value problem
\begin{equation}\label{d2}
\left\{\begin{array}{ll}
v_{1}-i\mu u_{1}=f_{1}&\text{in }\Omega_{1}
\\
v_{2}-i\mu u_{2}=f_{2}&\text{in }\Omega_{2}
\\
-\Delta(c_{1}^{2}\Delta u_{1}+a\Delta v_{1})-i\mu v_{1}=g_{1}&\text{in }\Omega_{1}
\\
-c_{2}^{2}\Delta^{2}u_{2}-i\mu v_{2}=g_{2}&\text{in }\Omega_{2}
\\
u_{1}=u_{2},\quad\pn u_{1}=\pn u_{2}&\text{on }S
\\
c_{1}\Delta u_{1}=c_{2}\Delta u_{2},\quad c_{1}\pn \Delta u_{1}=c_{2}\pn\Delta u_{2}&\text{on }S
\\
u_{2}=0,\quad\Delta u_{2}=0&\text{on }\Gamma.
\end{array}\right.
\end{equation}
Then the solution $(u_{1},u_{2},v_{1},v_{2})$ of~\eqref{d2} satisfies
\begin{equation}\label{d3}
\left\{\begin{array}{ll}
v_{1}=i\mu u_{1}+f_{1}&\text{in }\Omega_{1}
\\
v_{2}=i\mu u_{2}+f_{2}&\text{in }\Omega_{2}
\\
\mu^{2} u_{1}-\Delta(c_{1}^{2}\Delta u_{1}+a\Delta v_{1})=g_{1}+i\mu f_{1}&\text{in }\Omega_{1}
\\
\mu^{2}u_{2}-c_{2}^{2}\Delta^{2}u_{2}=g_{2}+i\mu f_{2}&\text{in }\Omega_{2}
\\
u_{1}=u_{2},\quad\pn u_{1}=\pn u_{2}&\text{on }S
\\
c_{1}\Delta u_{1}=c_{2}\Delta u_{2},\quad c_{1}\pn \Delta u_{1}=c_{2}\pn\Delta u_{2}&\text{on }S
\\
u_{2}=0,\quad\Delta u_{2}=0&\text{on }\Gamma.
\end{array}\right.
\end{equation}
This can be rewriten as follows
\begin{equation}\label{d4}
\left\{\begin{array}{ll}
v_{1}=i\mu u_{1}+f_{1}&\text{in }\Omega_{1}
\\
v_{2}=i\mu u_{2}+f_{2}&\text{in }\Omega_{2}
\\
\displaystyle(-\Delta-\frac{|\mu|}{c_{1}})(c_{1}\Delta u_{1}+\frac{a}{c_{1}}\Delta v_{1}-|\mu| u_{1})=\Phi_{1}=\frac{1}{c_{1}}g_{1}+i\frac{\mu}{c_{1}}f_{1}-a\frac{|\mu|}{c_{1}^{2}}\Delta v_{1}&\text{in }\Omega_{1}
\\
\\
\displaystyle(-\Delta-\frac{|\mu|}{c_{2}})(c_{2}\Delta u_{2}-|\mu|u_{2})=\Phi_{2}=\frac{1}{c_{2}}g_{2}+i\frac{\mu}{c_{2}}f_{2}&\text{in }\Omega_{2}
\\
u_{1}=u_{2},\quad\pn u_{1}=\pn u_{2}&\text{on }S
\\
c_{1}\Delta u_{1}=c_{2}\Delta u_{2},\quad c_{1}\pn \Delta u_{1}=c_{2}\pn\Delta u_{2}&\text{on }S
\\
u_{2}=0,\quad\Delta u_{2}=0&\text{on }\Gamma.
\end{array}\right.
\end{equation}
We set now 
\begin{equation}\label{d5}
w_{1}=c_{1}\Delta u_{1}-|\mu|u_{1}+\frac{a}{c_{1}}\Delta v_{1}\quad\text{ and }\quad w_{2}=c_{2}\Delta u_{2}-|\mu|u_{2},
\end{equation}
then it easy to show that $w_{1}$ and $w_{2}$ satisfy the following simple transmission problem
\begin{equation}\label{d6}
\left\{\begin{array}{ll}
\displaystyle-\Delta w_{1}-\frac{|\mu|}{c_{1}}w_{1}=\Phi_{1}&\text{in }\Omega_{1}
\\
\\
\displaystyle-\Delta w_{2}-\frac{|\mu|}{c_{2}}w_{2}=\Phi_{2}&\text{in }\Omega_{2}
\\
w_{1}=w_{2},\quad\pn w_{1}=\pn w_{2}&\text{on }S
\\
w_{2}=0&\text{on }\Gamma.
\end{array}\right.
\end{equation}
We set also $B_{4r}$ a ball of raduis $4r>0$, such that $a(x)>0$ in $B_{4r}\subset\omega$ and we recall the notation gived in the end of the previous section $\widetilde{\Omega}_{1}=\Omega_{1}\backslash \overline{B}_{r}$. The most important ingredient of the proof of the resolvent estimate is the following lemma which is essentially a consequence of the Carleman estimate.
\begin{lem}
There exist a constant $C>0$ such that for any $(u_{1},u_{2},v_{1},v_{2})\in\mathcal{D}(\mathcal{A})$ solution of~\eqref{d2} the following result holds
\begin{equation}\label{d7}
\begin{split}
\|\Delta u_{1}\|_{L^{2}(\Omega_{1})}^{2}+\|\Delta u_{2}\|_{L^{2}(\Omega_{2})}^{2}+\|v_{1}\|_{L^{2}(\Omega_{1})}^{2}+\|v_{2}\|_{L^{2}(\Omega_{2})}^{2}\leq C\e^{C|\mu|}\bigg(\|\Delta f_{1}\|_{L^{2}(\Omega_{1})}^{2}
\\
+\|\Delta f_{2}\|_{L^{2}(\Omega_{2})}^{2}+\|g_{1}\|_{L^{2}(\Omega_{1})}^{2}+\|g_{2}\|_{L^{2}(\Omega_{2})}^{2}+\int_{\Omega_{1}}a|\Delta v_{1}|^{2}\,\ud x+\int_{B_{4r}}|u_{1}|^{2}\,\ud x\bigg),
\end{split}
\end{equation}
for all $\mu\in\R$ large enough.
\end{lem}
\begin{pr}
We introduce the cutt-off function $\chi\in\Ci(\Omega_{1})$ by setting
$$
\chi(x)=\left\{\begin{array}{ll}
1&\text{in }B_{3r}^{c}
\\
0&\text{in }B_{2r}
\end{array}\right.
$$
Next, denote $\tilde{w}_{1}=\chi w_{1}$. And by~\eqref{d6}, one sees that
\begin{equation}\label{d8}
-\Delta\tilde{w}_{1}-\frac{|\mu|}{c_{1}}\tilde{w}_{1}=\widetilde{\Phi}_{1}=\chi\Phi_{1}-[\Delta,\chi]w_{1}.
\end{equation}

Now keeping the same notations as the previous section and let $\foo$, $\fot$, $\fto$ and $\ftt$ four weight functions that satisfies the conclusion of the section~\ref{b1}. Let $\chi_{1,1}$, $\chi_{1,2}$, $\chi_{2,1}$ and $\chi_{2,2}$ four cut-off functions equal to one respectively in $\displaystyle\left(\bigcup_{j=1}^{N_{11}}B(x_{1j}^{1},2\epsilon)\right)^{c}$, $\displaystyle\left(\bigcup_{j=1}^{N_{12}}B(x_{1j}^{2},2\epsilon)\right)^{c}$, $\displaystyle\left(\bigcup_{j=1}^{N_{21}}B(x_{2j}^{1},2\epsilon)\right)^{c}$ and $\displaystyle\left(\bigcup_{j=1}^{N_{22}}B(x_{2j}^{2},2\epsilon)\right)^{c}$ and supported respectively in $\displaystyle\left(\bigcup_{j=1}^{N_{11}}B(x_{1j}^{1},\epsilon)\right)^{c}$, $\displaystyle\left(\bigcup_{j=1}^{N_{12}}B(x_{1j}^{2},\epsilon)\right)^{c}$, $\displaystyle\left(\bigcup_{j=1}^{N_{21}}B(x_{2j}^{1},\epsilon)\right)^{c}$ and $\displaystyle\left(\bigcup_{j=1}^{N_{22}}B(x_{2j}^{2},\epsilon)\right)^{c}$ (in order to eliminate the critical points of the phases functions $\foo$, $\fot$, $\fto$ and $\ftt$ (See Figure~\ref{fig2})). We set now $w_{1,1}=\chi_{1,1}\tilde{w}_{1}$, $w_{1,2}=\chi_{1,2}\tilde{w}_{1}$, $w_{2,1}=\chi_{2,1}w_{2}$ and $w_{2,2}=\chi_{2,2}w_{2}$. Then from the system~\eqref{d6} for $k=1,2$ we obtain
\begin{equation}\label{d9}
\left\{\begin{array}{ll}
\displaystyle-\Delta w_{1,k}-\frac{|\mu|}{c_{1}}w_{1,k}=\Psi_{1,k}&\text{in }\Omega_{1}
\\
\\
\displaystyle-\Delta w_{2,k}-\frac{|\mu|}{c_{2}}w_{2,k}=\Psi_{2,k}&\text{in }\Omega_{2}
\\
w_{1,k}=w_{2,k},\quad\pn w_{1,k}=\pn w_{2,k}&\text{on }S
\\
w_{2,k}=0&\text{on }\Gamma,
\end{array}\right.
\end{equation}
where
\begin{equation}\label{d10}
\left\{\begin{array}{l}
\Psi_{1,k}=\chi_{1,k}\widetilde{\Phi}_{1}-[\Delta,\chi_{1,k}]\tilde{w}_{1}
\\
\Psi_{2,k}=\chi_{2,k}\Phi_{2}-[\Delta,\chi_{2,k}]w_{2}.
\end{array}\right.
\end{equation}
Applying now the Carleman estimate gived in the previous section (Theorem~\ref{b3}) to the system~\eqref{d9} for $\displaystyle h=\frac{1}{|\mu|}$ then for $k=1,2$ we obtain
\begin{equation*}
\begin{split}
h\|\e^{\vfi_{1,k}/h}w_{1,k}\|_{L^{2}(U_{1,k})}^{2}+h^{3}\|\e^{\vfi_{1,k}/h}\nabla w_{1,k}\|_{L^{2}(U_{1,k})}^{2}+h\|\e^{\vfi_{2,k}/h}w_{2,k}\|_{L^{2}(U_{2,k})}^{2}+
\\
h^{3}\|\e^{\vfi_{2,k}/h}\nabla w_{2,k}\|_{L^{2}(U_{2,k})}^{2}\leq Ch^{4}(\|\e^{\vfi_{1,k}/h}\Psi_{1,k}\|_{L^{2}(U_{1,k})}^{2}+\|\e^{\vfi_{2,k}/h}\Psi_{2,k}\|_{L^{2}(U_{2,k})}^{2}).
\end{split}
\end{equation*}
Relations~\eqref{d8} and~\eqref{d10} yields
\begin{equation}\label{d11}
\begin{split}
h\|\e^{\vfi_{1,k}/h}w_{1,k}\|_{L^{2}(U_{1,k})}^{2}+h^{3}\|\e^{\vfi_{1,k}/h}\nabla w_{1,k}\|_{L^{2}(U_{1,k})}^{2}+h\|\e^{\vfi_{2,k}/h}w_{2,k}\|_{L^{2}(U_{2,k})}^{2}+
\\
h^{3}\|\e^{\vfi_{2,k}/h}\nabla w_{2,k}\|_{L^{2}(U_{2,k})}^{2}\leq Ch^{4}(\|\e^{\vfi_{1,k}/h}\Phi_{1}\|_{L^{2}(U_{1,k})}^{2}+\|\e^{\vfi_{2,k}/h}\Phi_{2}\|_{L^{2}(U_{2,k})}^{2}+
\\
\|\e^{\vfi_{1,k}/h}[\Delta,\chi]w_{1}\|_{L^{2}(U_{1,k})}^{2}+\|\e^{\vfi_{1,k}/h}[\Delta,\chi_{1,k}]\tilde{w}_{1}\|_{L^{2}(U_{1,k})}^{2}+\|\e^{\vfi_{2,k}/h}[\Delta,\chi_{2,k}]w_{2}\|_{L^{2}(U_{2,k})}^{2}).
\end{split}
\end{equation}
We addition the two last estimates  for $k=1,2$ and using the properties of phases $\vfi_{1,k}<\vfi_{1,\sigma(k)}$ in $\displaystyle\left(\bigcup_{j=1}^{N_{1k}}B(x_{1k}^{j},2\epsilon)\right)$ and $\vfi_{2,k}<\vfi_{2,\sigma(k)}$ in $\displaystyle\left(\bigcup_{j=1}^{N_{2k}}B(x_{2k}^{j},2\epsilon)\right)$ then we can absorb the terms $[\Delta,\chi_{1,k}]\tilde{w}_{1}$ and $[\Delta,\chi_{2,k}]w_{2}$ at the right hand side of~\eqref{d11} into the left hand side for $h>0$ small. More precisly we obtain
\begin{equation*}
\begin{split}
h\int_{\widetilde{\Omega}_{1}}(\e^{2\foo/h}+\e^{2\fot/h})|\tilde{w}_{1}|^{2}\,\ud x+h\int_{\Omega_{2}}(\e^{2\fto/h}+\e^{2\ftt/h})|w_{2}|^{2}\,\ud x\leq
\\
Ch^{4}\bigg(\int_{\Omega_{1}}(\e^{2\foo/h}+\e^{2\fot/h})|\Phi_{1}|^{2}\,\ud x+\int_{\Omega_{2}}(\e^{2\fto/h}+\e^{2\ftt/h})|\Phi_{2}|^{2}\,\ud x
\\
+\int_{\widetilde{\Omega}_{1}}(\e^{2\foo/h}+\e^{2\fot/h})|[\Delta,\chi]w_{1}|^{2}\,\ud x\bigg).
\end{split}
\end{equation*}
Consequently, by using that $\Omega_{1}=\widetilde{\Omega}_{1}\cup B_{2r}$ and the expressions of $\Phi_{1}$ and $\Phi_{2}$ in~\eqref{d4} we see that
\begin{equation}\label{d12}
\begin{split}
\int_{\Omega_{1}}|w_{1}|^{2}\,\ud x+\int_{\Omega_{2}}|w_{2}|^{2}\,\ud x\leq C\e^{C/h}\bigg(\int_{\Omega_{1}}|f_{1}|^{2}\,\ud x+\int_{\Omega_{1}}|g_{1}|^{2}\,\ud x+\int_{\Omega_{2}}|f_{2}|^{2}\,\ud x
\\
+\int_{\Omega_{2}}|g_{2}|^{2}\,\ud x
+\int_{\Omega_{1}}a|\Delta v_{1}|^{2}\,\ud x+\int_{B_{2r}}|w_{1}|^{2}\,\ud x+\int_{\widetilde{\Omega}_{1}}|[\Delta,\chi]w_{1}|^{2}\,\ud x\bigg).
\end{split}
\end{equation}
To accomplish the proof of the lemma we estimate the two last terms in the right hand side of~\eqref{d12}. We set $\widetilde{\chi}$ a cutt-off function equal to $1$ in a neighborhood of $B_{3r}$ and supported in $B_{4r}$ then we have
$$
(-1+\Delta)(\widetilde{\chi}w_{1})=[\Delta,\widetilde{\chi}]w_{1}-\widetilde{\chi}w_{1}-\frac{|\mu|}{c_{1}}\widetilde{\chi}w_{1}-\widetilde{\chi}\Phi_{1},
$$
and hence by elliptic estimates (see~\cite{WRL}) we get
\begin{eqnarray}\label{d13}
\|w_{1}\|_{H^{1}(B_{3r})}^{2}\!\!\!\!&\leq&\!\!\!C(\|(-1+\Delta)(\widetilde{\chi}w_{1})\|_{H^{-1}(B_{4r})}^{2}+\|w_{1}\|_{L^{2}(B_{4r})}^{2})\nonumber
\\
&\leq&\!\!\!C(\|\Phi_{1}\|_{L^{2}(\Omega_{1})}^{2}+(1+|\mu|^{2})\|w_{1}\|_{L^{2}(B_{4r})}^{2})\nonumber
\\
&\leq&\!\!\!C\left(|\mu|^{2}\|f_{1}\|_{L^{2}(\Omega_{1})}^{2}+\|g_{1}\|_{L^{2}(\Omega_{1})}^{2}+(1+|\mu|^{2})\|w_{1}\|_{L^{2}(B_{4r})}^{2}+|\mu|^{2}\int_{\Omega_{1}}a|\Delta v_{1}|^{2}\,\ud x\right).
\end{eqnarray}
Since $\supp([\Delta,\chi])\subset B_{3r}$ we deduce from~\eqref{d5} and~\eqref{d13} that
\begin{equation}\label{d14}
\begin{split}
\int_{B_{2r}}|w_{1}|^{2}\,\ud x+\int_{\widetilde{\Omega}_{1}}|[\Delta,\chi]w_{1}|^{2}\,\ud x\leq C\|w_{1}\|_{H^{1}(B_{3r})}^{2}
\\
\leq C\left(|\mu|^{2}\|f_{1}\|_{L^{2}(\Omega_{1})}^{2}+\|g_{1}\|_{L^{2}(\Omega_{1})}^{2}+(1+|\mu|^{2})^{2}\|u_{1}\|_{L^{2}(B_{4r})}^{2}+|\mu|^{2}\int_{\Omega_{1}}a|\Delta v_{1}|^{2}\,\ud x\right).
\end{split}
\end{equation}
On other hand from~\eqref{d5} and the transmission conditions we see that
\begin{equation}\label{d15}
\begin{split}
\|w_{1}\|_{L^{2}(\Omega_{1})}^{2}+\|w_{2}\|_{L^{2}(\Omega_{2})}^{2}\geq \|c_{1}\Delta u_{1}-|\mu|u_{1}\|_{L^{2}(\Omega_{1})}^{2}+\|c_{2}\Delta u_{2}-|\mu|u_{2}\|_{L^{2}(\Omega_{2})}^{2}-C\int_{\Omega_{1}}\!\!\!a|\Delta v_{1}|^{2}\,\ud x
\\
\geq -C\int_{\Omega_{1}}\!\!\!a|\Delta v_{1}|^{2}\,\ud x+c_{1}^{2}\|\Delta u_{1}\|_{L^{2}(\Omega_{1})}^{2}+c_{2}^{2}\|\Delta u_{2}\|_{L^{2}(\Omega_{2})}^{2}+|\mu|^{2}(\|u_{1}\|_{L^{2}(\Omega_{1})}^{2}+\|u_{2}\|_{L^{2}(\Omega_{2})}^{2})
\\
+|\mu|(\|\nabla u_{1}\|_{L^{2}(\Omega_{1})}^{2}+\|\nabla u_{2}\|_{L^{2}(\Omega_{2})}^{2})\geq\|\Delta u_{1}\|_{L^{2}(\Omega_{1})}^{2}+\|\Delta u_{2}\|_{L^{2}(\Omega_{2})}^{2}-C\int_{\Omega_{1}}\!\!\!a|\Delta v_{1}|^{2}\,\ud x,
\end{split}
\end{equation}
and by the expression of $v_{1}$ and $v_{2}$ in~\eqref{d3} we obtain
\begin{equation}\label{d16}
\begin{split}
&\|v_{1}\|_{L^{2}(\Omega_{1})}^{2}\leq \|f_{1}\|_{L^{2}(\Omega_{1})}^{2}+|\mu|^{2}\|u_{1}\|_{L^{2}(\Omega_{1})}^{2}
\\
&\|v_{2}\|_{L^{2}(\Omega_{2})}^{2}\leq \|f_{2}\|_{L^{2}(\Omega_{2})}^{2}+|\mu|^{2}\|u_{2}\|_{L^{2}(\Omega_{2})}^{2}.
\end{split}
\end{equation}
Then by combining Proposition~\ref{7}, and estimates~\eqref{d12},~\eqref{d14},~\eqref{d15} and~\eqref{d16} we obtain the results.
\end{pr}

At this step we suppose now that the resolvent estimate~\eqref{2} is not true. Then there exist $K_{m}>0$, $\mu_{m}\in\R$ and a two families $(u_{1,m},u_{2,m},v_{1,m},v_{2,m})\in\mathcal{D}(\mathcal{A})$ and $(f_{1,m},f_{2,m},g_{1,m},g_{2,m})\in\mathcal{H}$, $m=1,2,\ldots$  such that
\begin{equation}\label{d17}
|\mu_{m}|\,\longrightarrow\,+\infty,\qquad K_{m}\,\longrightarrow\,+\infty,\qquad\|(u_{1,m},u_{2,m},v_{1,m},v_{2,m})\|_{\mathcal{H}}=1,
\end{equation}
and
\begin{equation}\label{d18}
\e^{K_{m}|\mu_{m}|}(\mathcal{A}-i\mu_{m})\left(\begin{array}{c}
u_{1,m}
\\
u_{2,m}
\\
v_{1,m}
\\
v_{2,m}
\end{array}\right)=\left(\begin{array}{c}
f_{1,m}
\\
f_{2,m}
\\
g_{1,m}
\\
g_{2,m}
\end{array}\right)\,\longrightarrow\,0\text{ in }\mathcal{H}.
\end{equation}
This imply that
\begin{eqnarray}
\e^{K_{m}|\mu_{m}|}(v_{1,m}-i\mu_{m}u_{1,m})=f_{1,m}\,\longrightarrow\,0\text{ in }  H^{2}(\Omega_{1}),\label{d19}
\\
\e^{K_{m}|\mu_{m}|}(v_{2,m}-i\mu_{m}u_{2,m})=f_{2,m}\,\longrightarrow\,0\text{ in }  H^{2}(\Omega_{2}),
\\
\e^{K_{m}|\mu_{m}|}(-\Delta(c_{1}^{2}\Delta u_{1,m}+a\Delta v_{1,m})-i\mu_{m}v_{1,m})=g_{1,m}\,\longrightarrow\,0\text{ in }  L^{2}(\Omega_{1}),\label{d20}
\\
\e^{K_{m}|\mu_{m}|}(-c_{2}^{2}\Delta^{2}u_{2,m}-i\mu_{m}v_{2,m})=g_{2,m}\,\longrightarrow\,0\text{ in } L^{2}(\Omega_{2}).\label{d26}
\end{eqnarray}
From~\eqref{d17} and~\eqref{d18}, we get
\begin{equation}\label{d21}
\re\left\langle\left(\begin{array}{c}
f_{1,m}
\\
f_{2,m}
\\
g_{1,m}
\\
g_{2,m}
\end{array}\right),\left(\begin{array}{c}
u_{1,m}
\\
u_{2,m}
\\
v_{1,m}
\\
v_{2,m}
\end{array}\right)\right\rangle_{\mathcal{H}}=-\e^{K_{m}|\mu_{m}|}\int_{\Omega_{1}}a|\Delta v_{1,m}|^{2}\,\ud x\,\longrightarrow\,0.
\end{equation}
Then by~\eqref{d19} and~\eqref{d21}, we obtain
\begin{equation}\label{d22}
|\mu_{m}|^{2}\e^{\frac{K_{m}}{2}|\mu_{m}|}\int_{\omega}|\Delta u_{1,m}|^{2}\,\ud x\,\longrightarrow\,0.
\end{equation}
Hence from~\eqref{d21} and~\eqref{d22} we obtain
\begin{equation}\label{d23}
\e^{\frac{K_{m}}{2}|\mu_{m}|}\left(\int_{\omega}|\Delta u_{1,m}|^{2}\,\ud x+\int_{\omega}|\Delta v_{1,m}|^{2}\,\ud x\right)\,\longrightarrow\,0.
\end{equation}
And by~\eqref{d19} we have
\begin{equation}\label{d24}
\frac{1}{|\mu_{m}|^{2}}\|\Delta(\psi.v_{1,m})\|_{L^{2}(\Omega_{1})}^{2}=O(1),\quad\forall\,\psi\in\Ci(\Omega_{1}).
\end{equation}
Then by multiplying~\eqref{d20} by $\mu_{m}^{-1}\psi.\overline{v}_{1,m}$ where $\psi\in\Ci(\Omega_{1})$ and $\supp(\psi)\subset\omega$ we obtain by ~\eqref{d23} and~\eqref{d24} that
\begin{equation*}
\e^{\frac{K_{m}}{4}|\mu_{m}|}\int_{\omega}|v_{1,m}|^{2}\psi\,\ud x\,\longrightarrow\,0.
\end{equation*}
In particular we obtain that
\begin{equation*}
\e^{\frac{K_{m}}{4}|\mu_{m}|}\int_{B_{4r}}|v_{1,m}|^{2}\,\ud x\,\longrightarrow\,0.
\end{equation*}
Then also we get by~\eqref{d19} that
\begin{equation}\label{d25}
\e^{\frac{K_{m}}{4}|\mu_{m}|}\int_{B_{4r}}|u_{1,m}|^{2}\,\ud x\,\longrightarrow\,0.
\end{equation}
Now by applying inequality~\eqref{d7} to the system~\eqref{d19}-\eqref{d26} it follows that
\begin{equation}\label{d27}
\begin{split}
\|\Delta u_{1,m}\|_{L^{2}(\Omega_{1})}^{2}+\|\Delta u_{2,m}\|_{L^{2}(\Omega_{2})}^{2}+\|v_{1,m}\|_{L^{2}(\Omega_{1})}^{2}+\|v_{2,m}\|_{L^{2}(\Omega_{2})}^{2}\leq
\\
C\e^{C|\mu_{m}|}\bigg(\e^{-2K_{m}|\mu_{m}|}\Big(\|\Delta f_{1,m}\|_{L^{2}(\Omega_{1})}^{2}+\|\Delta f_{2,m}\|_{L^{2}(\Omega_{2})}^{2}+\|g_{1,m}\|_{L^{2}(\Omega_{1})}^{2}+\|g_{2,m}\|_{L^{2}(\Omega_{2})}^{2}\Big)
\\
+\e^{-\frac{K_{m}}{4}|\mu_{m}|}\left(\int_{\Omega_{1}}a|\Delta v_{1,m}|^{2}\,\ud x+\int_{B_{4r}}|u_{1,m}|^{2}\,\ud x\right)\e^{\frac{K_{m}}{4}|\mu_{m}|}\bigg).
\end{split}
\end{equation}
While the right hand side of~\eqref{d27} go to zero as $m\,\longrightarrow\,+\infty$ by~\eqref{d17}-\eqref{d18} and estimates~\eqref{d21} and~\eqref{d25}, then we obtain a contradiction with~\eqref{d17}. And this conclude the proof of the resolvent estimate.
\nocite{*}
\bibliographystyle{alpha}
\bibliography{BibKVMultyD}

\begin{thebibliography}{CFNS91}

\bibitem[Alb00]{A}
P.~Albano.
\newblock Carleman estimates for the {E}uler-{B}ernoulli plate operator.
\newblock {\em Electronic journal of differential equations}, pages 1--13,
  2000.

\bibitem[AN10]{AN}
K.~Ammari and S.~Nicaise.
\newblock Stabilization of a transmission wave/plate equation.
\newblock {\em Journal of Differential Equations}, 249:707--727, 2010.

\bibitem[AV09]{AV}
K.~Ammari and G.~Vodev.
\newblock Boundary stabilization of the transmission problem for the
  {B}ernoulli-{E}uler plate equation.
\newblock {\em CUBO a mathematical journal}, 11:39--49, 2009.

\bibitem[BD08]{BD}
C.J.K Batty and T.~Duyckaerts.
\newblock Non-uniform stability for bounded semi-groups on {B}anach spaces.
\newblock {\em Journal of Evolution Equation}, pages 765--780, 2008.

\bibitem[Bel03]{B}
M.~Bellassoued.
\newblock Carleman estimates and distribution of resonnances for the
  transparent obstace and application to the stabilization.
\newblock {\em Asymptotic Anal.}, 35:257--279, 2003.

\bibitem[Bur98]{Bur}
N.~Burq.
\newblock D\'ecroissance de l'\'energie locale de l'\'equation des ondes pour
  le probl\`eme ext\'erieur et absence de r\'esonnance au voisinage du r\'eel.
\newblock {\em Acta Math.}, 180:1--29, 1998.

\bibitem[CFNS91]{CFNS}
G.~Chen, S.~A. Fulling, F.~J. Narcowich, and S.~Sun.
\newblock Exponential decay of energy of evolution equations with locally
  distributed damping.
\newblock {\em SIAM J. Appl. Math.}, 51(1):266--301, 1991.

\bibitem[CLL98]{CLL}
S.~Chen, K.~Liu, and Z.~Liu.
\newblock Spectrum and stability for elastic systems with global or local
  {K}elvin-{V}oigt damping.
\newblock {\em SIAM J. APPL. MATH.}, 59(2):651--668, 1998.

\bibitem[EMT04]{YVA}
Y.~Eidelman, V.~Milman, and A.~Tsolomitis.
\newblock {\em Functional Analysis An Introduction}.
\newblock American Mathematical Society, 2004.

\bibitem[Fat11]{I}
I.K. Fathallah.
\newblock Logarithmic decay of the energy for an hyperbolic-parabolic coupled
  system.
\newblock {\em ESAIM-control Optimisation and Calculus of Variations},
  17:801--835, 2011.

\bibitem[Lau12]{L}
F.~Laudenbach.
\newblock A proof of {R}eidemeister-{S}inger's theorem by {C}erf's methods.
\newblock {\em ArXiv(math.GT):1202.1130}, 2012.

\bibitem[LL98]{LL}
K.~Liu and Z.~Liu.
\newblock Exponential decay of the energy of the {E}uler-{B}ernoulli beam with
  locally distributed {K}elvin-{V}oigt damping.
\newblock {\em SIAM J. Control Optim.}, 36(3):1086--1098, 1998.

\bibitem[LL02]{LL2}
K.~Liu and Z.~Liu.
\newblock Exponential decay of energy of vibrating strings with local
  viscoelasticity.
\newblock {\em Z. Angew Math. Phys.}, 53:265--280, 2002.

\bibitem[LR95]{LR2}
G.~Lebeau and L.~Robbiano.
\newblock Contr\^ole exacte de l'\'equation de la chaleur.
\newblock {\em Comm. Partial Differential Equations}, 20:335--356, 1995.

\bibitem[LR97]{LR}
G.~Lebeau and L.~Robbiano.
\newblock Stabilisation de l'\'equation des ondes par le bord.
\newblock {\em Duke mathematical journal}, 86:465--491, december 1997.

\bibitem[Mil65]{Mi}
J.~Milnor.
\newblock {\em Lectures on the h-{C}obordism {T}heorem}.
\newblock Princeton Univ. Press, Princeton, NJ, 1965.

\bibitem[RR10]{RR}
J.~Le Rousseau and L.~Robbiano.
\newblock Carleman estimate for elliptic operators with coefficients with jumps
  at an interface in arbitrary dimension and application to the null
  controllability of linear parabolic equations.
\newblock {\em Arch. Rational Mech. Anal.}, 195:953--990, 2010.

\bibitem[TW09]{TW}
M.~Tucsnak and G.~Weiss.
\newblock {\em Observation and control for operator semigroups}.
\newblock Birkh\"auser Verlag AG, 2009.

\bibitem[WRL95]{WRL}
J.~T. Wolka, B.~Rowley, and B.~Lawruk.
\newblock {\em Boundary value problems for elliptic system}.
\newblock Cambridge University Press, Cambridge, 1995.

\bibitem[Zua90]{Z}
E.~Zuazua.
\newblock Exponential decay for the semilinear wave equation with localized
  damping.
\newblock {\em Comm. Part. Diff. Eq.}, 15:205--235, 1990.

\end{thebibliography}
\addcontentsline{toc}{section}{References}
\end{document}